\documentclass{amsart}
\usepackage{amssymb}
\usepackage{amsmath}
\usepackage{amstext}
\usepackage{amsfonts}
\usepackage{amsthm}
\usepackage{mathrsfs}
\usepackage{centernot}
\usepackage{enumitem}
\usepackage{amsaddr}

\newtheorem{thm}{Theorem}[section]
\newtheorem{thmx}{Theorem}[section]

\theoremstyle{definition}
\newtheorem{defn}[thm]{Definition}
\newtheorem{defnx}[thmx]{Definition}
\theoremstyle{remark}
\newtheorem{rem}[thm]{Remark}
\theoremstyle{example}
\newtheorem{exam}[thm]{Example}
\numberwithin{equation}{section}
\makeatother

\begin{document}
\title{Strong boundedness, strong convergence and generalized variation}
\author{M. Avdispahi\'{c} \MakeLowercase{and} Z. \v{S}abanac}
\address{\textnormal{Department of Mathematics, University of Sarajevo,
Zmaja od Bosne 33-35 Sarajevo, Bosnia and Herzegovina}\\
\textnormal{e-mails: mavdispa@pmf.unsa.ba, zsabanac@pmf.unsa.ba}}
\maketitle

\dedicatory{\textit{\hfill To the memory of Naza Tanovi\'{c}-Miller}}

\begin{abstract}
A trigonometric series strongly bounded at two points and with coefficients forming a log-quasidecreasing sequence is necessarily the Fourier series of a function belonging to all $L^{p}$ spaces, $1\leq p < \infty$.  We obtain new results on strong convergence of Fourier series for functions of generalized bounded variation.
\end{abstract}

\section{Introduction}
Extending Hardy-Littlewood's concept of strong $(C,1)$ summability to Ces%
\`{a}ro methods $(C,\alpha )$ of order $\alpha \geq 0$, Hyslop \cite{Hy}
arrived at his notion of strong convergence. Subsequently, this was
successfully applied to the study of trigonometric series in several papers
written by N. Tanovi\'{c}-Miller and her co-workers \cite{NTM1, NTM2, NTM3,
STM, STM1}. Strong convergence of trigonometric series attracts attention
because of its position between ordinary and absolute convergence \cite{AV2,
NTM1, NTM2}.

Interesting results about the global behaviour of a series deduced from its
behaviour at one or two points were initially related to absolute
convergence and obtained by O. Szasz \cite{SZ} and R. Pippert \cite{PI}. The
assumption on the coefficients of a series was that their magnitudes form an
almost decreasing sequence. The analogues are valid in the case of strong
convergence \cite{AV1}.

We introduce new notions of strong boundedness (in Hyslop's sense) and
logarithmic quasimonotonicity. We prove that if one requests only strong
boundedness of a trigonometric series at two points but imposes logarithmic
quasimonotonicity on the magnitudes of its coefficients, then the respective
trigonometric series is the Fourier series of a function belonging to all $%
L^{p}$ spaces, $1\leq p<\infty $.

In the area of strong convergence, our attention is turned to Fourier series
of regulated functions, i.e., functions belonging to various classes of
generalized bounded variation.

\section{Banach spaces of strongly bounded sequences }

\begin{defnx}
A sequence of numbers $\left\{ d_{n}\right\} $ is strongly $\left(
C,1\right) $ summable to a limit $d$ with index $\lambda >0$ ($\lambda -$%
strongly $\left( C,1\right) $ summable to $d$), and we write $%
d_{n}\rightarrow d \ \left[ C_{1}\right] _{\lambda }$, if%
\begin{equation*}
\sum_{k=1}^{n}\left\vert d_{k}-d\right\vert ^{\lambda }=o\left( n\right)
\text{ as }n\rightarrow \infty \text{.}
\end{equation*}
\end{defnx}

\begin{defnx}
A sequence of numbers $\left\{ d_{n}\right\} $ is strongly convergent to a
limit $d$ with index $\lambda >0$ ($\lambda -$strongly convergent to $d$),
and we write $d_{n}\rightarrow d \ \left[ I\right] _{\lambda }$, if

\begin{description}
\item[1)] $d_{n}\rightarrow d$ as $n\rightarrow \infty $,

\item[2)] $\sum_{k=1}^{n}k^{\lambda }\left\vert d_{k}-d_{k-1}\right\vert
^{\lambda }=o\left( n\right) $ as $n\rightarrow \infty $, i.e., $k\left(
d_{k}-d_{k-1}\right) \rightarrow 0 \ \left[ C_{1}\right] _{\lambda }$.
\end{description}

If $\lambda =1$, we simply denote it by $\left[ I\right] $.
\end{defnx}

\begin{defn}
A sequence of numbers $\left\{ d_{n}\right\} $ is said to be strongly
bounded with index $\lambda >0$ ($\lambda -$strongly bounded), if

\begin{description}
\item[1)] $d_{n}=O\left( 1\right) $ as $n\rightarrow \infty $,

\item[2)] $\sum_{k=1}^{n}k^{\lambda }\left\vert d_{k}-d_{k-1}\right\vert
^{\lambda }=O\left( n\right) $ as $n\rightarrow \infty $.
\end{description}

If $\lambda =1$, we say that sequence $\left\{ d_{n}\right\} $ is
strongly bounded.

The set of $\lambda -$strongly bounded sequences is denoted by $\mathscr{B}%
^{\lambda }$.
\end{defn}

\bigskip

It is obvious that every $\lambda -$strongly convergent sequence is $\lambda
-$strongly bounded. The converse does not hold as illustrated by the following examples.

\begin{exam}
Let a sequence $\left\{ d_{n}\right\} $ be given by $%
d_{n}=\sum_{k=1}^{n}\frac{\left( -1\right) ^{k-1}}{k}$ for $n\in
\mathbb{N}
$. This sequence is obviously convergent.Therefore, $d_{n}=O\left( 1\right) $
as $n\rightarrow \infty$. For any $\lambda >0$, we have%
\begin{equation*}
\sum_{k=1}^{n}k^{\lambda }\left\vert d_{k}-d_{k-1}\right\vert ^{\lambda
}=\sum_{k=1}^{n}k^{\lambda }\left\vert \frac{\left( -1\right) ^{k-1}}{k}%
\right\vert ^{\lambda }=\sum_{k=1}^{n}1=n \text{.}
\end{equation*}
Hence, $\left\{ d_{n}\right\} $ is  $\lambda -$strongly bounded but not $\lambda -$strongly convergent.
\end{exam}

\begin{exam}
Consider a sequence $\left\{ d_{n}\right\} $ such that
\begin{equation*}
d_{n}=\left\{
\begin{array}{cl}
1\text{,} & \text{ if }n=2^{k}\text{, }k\in
\mathbb{N}
\text{,} \\
0\text{,} & \text{if }n\neq 2^{k}\text{, }k\in
\mathbb{N}
\text{.}%
\end{array}%
\right.
\end{equation*}%
This sequence is bounded but it is not convergent since it has two partial
limits, $0$ and $1$. Therefore, it is not $\lambda -$strongly convergent.
Now, we have $\left\vert d_{n}-d_{n-1}\right\vert =1$  if $n=2^{k}$ or $n=2^{k}+1$. Otherwise, $d_{n}-d_{n-1}=0$. Let $0<\lambda \leq1$. Then%
\begin{equation*}
\sum_{k=1}^{n}k^{\lambda }\left\vert d_{k}-d_{k-1}\right\vert ^{\lambda
}= \sum_{j=0}^{\lfloor \log_{2} n \rfloor}[2^{j\lambda}+(2^{j}+1)^{\lambda}]=O\left( n^{\lambda }\right)\text{.}
\end{equation*}
Therefore, this sequence is $\lambda -$strongly bounded for $0<\lambda \leq1$.
\end{exam}

Every $\mathscr{B}^{\lambda }$ is a linear space. The next theorem introduces a norm in $\mathscr{B}^{\lambda }$ that turns $\mathscr{B}^{\lambda }$ into a Banach space.

\begin{thm}
\leavevmode

\begin{itemize}
\item[i)] $\mathscr{B}^{\mu }\supseteq \mathscr{B}^{\lambda }$ for $0<\mu
<\lambda $.

\item[ii)] For $d=\left\{ d_{n}\right\} _{n=1}^{\infty }\in \mathscr{B}%
^{\lambda }$, $\lambda \geq 1$, let%
\begin{equation*}
\left\Vert d\right\Vert _{\mathscr{B}^{\lambda }}=\underset{n}{\sup }%
\left\vert d_{n}\right\vert +\underset{n}{\sup }\left( \frac{1}{n}%
\sum_{k=1}^{n}k^{\lambda }\left\vert d_{k}-d_{k-1}\right\vert ^{\lambda
}\right) ^{\frac{1}{\lambda }}\text{.}  \label{eq6}
\end{equation*}

$\left\Vert \cdot\right\Vert _{\mathscr{B}^{\lambda }}$ is a norm on $\mathscr{B}^{\lambda }$, $\lambda
\geq 1$.

\item[iii)] $\mathscr{B}^{\lambda }$, $\lambda \geq 1$, is a Banach space
under the norm given in $ii)$.
\end{itemize}

\begin{proof}
\mbox{}\\*[0pt]
$i)$ This is an immediate consequence of H\"{o}lder's inequality%
\begin{equation*}
\sum_{k=1}^{n}k^{\mu }\left\vert d_{k}-d_{k-1}\right\vert ^{\mu} \leq \left(\sum_{k=1}^{n}k^{\lambda }\left\vert d_{k}-d_{k-1}\right\vert ^{\lambda
}\right) ^{\frac{\mu}{\lambda }} \left(\sum_{k=1}^{n}1\right)^{1-\frac{\mu}{\lambda }}=O\left(n\right)
\end{equation*}%
for $0<\mu <\lambda $. \mbox{}\\*[%
0pt]
$ii)$ It is obvious that $\left\Vert d\right\Vert _{\mathscr{B}%
^{\lambda }}\geq 0$ and that the equality holds if and only if $d=\left\{
0\right\} _{n=1}^{\infty }$. If $\alpha $ is an arbitrary complex number
and $\alpha d:=\left\{ \alpha d_{n}\right\} _{n=1}^{\infty }$,
then
\begin{gather*}
\left\Vert \alpha d\right\Vert _{\mathscr{B}^{\lambda }}=\underset{n}{\sup }%
\left\vert \alpha d_{n}\right\vert +\underset{n}{\sup }\left( \frac{1}{n}%
\sum_{k=1}^{n}k^{\lambda }\left\vert \alpha d_{k}-\alpha d_{k-1}\right\vert
^{\lambda }\right) ^{\frac{1}{\lambda }} \\
=\left\vert \alpha \right\vert \left( \underset{n}{\sup }\left\vert
d_{k}\right\vert +\underset{n}{\sup }\left( \frac{1}{n}\sum_{k=1}^{n}k^{%
\lambda }\left\vert d_{k}-d_{k-1}\right\vert ^{\lambda }\right) ^{\frac{1}{%
\lambda }}\right) =\left\vert \alpha \right\vert \left\Vert d\right\Vert _{%
\mathscr{B}^{\lambda }}\text{.}
\end{gather*}%
If $d^{(1)}=\left\{ d_{n}^{(1)}\right\} _{n=1}^{\infty }$ and $%
d^{(2)}=\left\{ d_{n}^{(2)}\right\} _{n=1}^{\infty }$ are two $\lambda -$%
strongly bounded sequences, $\lambda \geq 1$, and $%
d^{(1)}+d^{(2)}:=\left\{ d_{n}^{(1)}+d_{n}^{(2)}\right\} _{n=1}^{\infty }$ ,
then by Minkowski's inequality we get%
\begin{gather*}
\left\vert d_{n}^{(1)}+d_{n}^{(2)}\right\vert +\left( \frac{1}{n}%
\sum_{k=1}^{n}k^{\lambda }\left\vert \left( d_{k}^{(1)}+d_{k}^{(2)}\right)
-\left( d_{k-1}^{(1)}+d_{k-1}^{(2)}\right) \right\vert ^{\lambda }\right) ^{%
\frac{1}{\lambda }} \\
=\left\vert d_{n}^{(1)}+d_{n}^{(2)}\right\vert +\left( \frac{1}{n}%
\sum_{k=1}^{n}k^{\lambda }\left\vert \left( d_{k}^{(1)}-d_{k-1}^{(1)}\right)
+\left( d_{k}^{(2)}-d_{k-1}^{(2)}\right) \right\vert ^{\lambda }\right) ^{%
\frac{1}{\lambda }} \\
\leq \left\vert d_{n}^{(1)}\right\vert +\left\vert d_{n}^{(2)}\right\vert
+\left( \frac{1}{n}\sum_{k=1}^{n}k^{\lambda }\left\vert
d_{k}^{(1)}-d_{k-1}^{(1)}\right\vert ^{\lambda }\right) ^{\frac{1}{\lambda }%
}+\left( \frac{1}{n}\sum_{k=1}^{n}k^{\lambda }\left\vert
d_{k}^{(2)}-d_{k-1}^{(2)}\right\vert ^{\lambda }\right) ^{\frac{1}{\lambda }}\text{.}
\end{gather*}%
Hence,%
\begin{gather*}
\left\Vert d^{(1)}+d^{(2)}\right\Vert _{\mathscr{B}^{\lambda }}\\
\leq \underset{n}{\sup }\left\vert d_{n}^{(1)}\right\vert +\underset{n}{\sup
}\left( \frac{1}{n}\sum_{k=1}^{n}k^{\lambda }\left\vert
d_{k}^{(1)}-d_{k-1}^{(1)}\right\vert ^{\lambda }\right) ^{\frac{1}{\lambda }%
}+\underset{n}{\sup }\left\vert d_{n}^{(2)}\right\vert +\underset{n}{\sup }%
\left( \frac{1}{n}\sum_{k=1}^{n}k^{\lambda }\left\vert
d_{k}^{(2)}-d_{k-1}^{(2)}\right\vert ^{\lambda }\right) ^{\frac{1}{\lambda }}
\\
=\left\Vert d^{(1)}\right\Vert _{\mathscr{B}^{\lambda }}+\left\Vert
d^{(2)}\right\Vert _{\mathscr{B}^{\lambda }}\text{.}
\end{gather*}%
Thus, $\mathscr{B}^{\lambda }$, $\lambda \geq
1$, is a normed linear space.\mbox{}\\*[0pt]
$iii)$ It remains to check that $\mathscr{B}^{\lambda }$, $\lambda \geq 1$,
is complete. Let $d^{(1)}$, $d^{(2)}$, \ldots\ be a Cauchy sequence in $%
\mathscr{B}^{\lambda }$. Now,
\begin{equation*}
\left( \forall \varepsilon >0\right) \left( \exists n_{0}\in
\mathbb{N}
\right) \left( \forall m>n\geq n_{0}\right) \left\Vert
d^{(m)}-d^{(n)}\right\Vert _{\mathscr{B}^{\lambda }}<\frac{\varepsilon}{3} \text{.}
\end{equation*}%
Note that
\begin{equation*}
\left\Vert
d^{(m)}-d^{(n)}\right\Vert _{l^{\infty }} \leq \left\Vert d^{(m)}-d^{(n)}\right\Vert _{\mathscr{B}^{\lambda }}\text{.}
\end{equation*}%
Thus, $d^{(n)}$ is a Cauchy sequence in $l^{\infty }$. Since $l^{\infty }$ is a Banach space, there exists $d=\left\{ d_{i}\right\} _{i=1}^{\infty }$ $\in l^{\infty }$
such that $\left\Vert
d^{(n)}-d\right\Vert _{l^{\infty }}\rightarrow 0 \ \ (n\rightarrow\infty)$. Hence, for $\varepsilon > 0$ chosen above,
\begin{equation}
\left( \exists n_{0}^{\ast }\in
\mathbb{N}
\right) \left( \forall n\geq n_{0}^{\ast }\right) \left\Vert
d^{(n)}-d\right\Vert _{l^{\infty }}<\frac{\varepsilon }{3}\text{.}
\label{eq7}
\end{equation}%
Moreover,%
\begin{equation}
\left( \forall k\in
\mathbb{N}
\right) \left( \exists n_{k}\in
\mathbb{N}
\right) \left( \forall n\geq n_{k}\right) \left\vert
d_{k}^{(n)}-d_{k}\right\vert <\frac{\varepsilon }{6\left( k+1\right) 2^{%
\frac{k+1}{\lambda }}}\text{.}  \label{eq8}
\end{equation}%
Let us show that $\left\{ d^{(n)}\right\} _{n=1}^{\infty }$ converges to $%
d=\left\{ d_{i}\right\} _{i=1}^{\infty }$ in $\mathscr{B}^{\lambda }$. Take an arbitrary $i \in \mathbb{N}$ and fix it. Put $n_{i}^{\ast }=\max
\left\{ n_{0}^{\ast },n_{1},n_{2},\ldots ,n_{i}\right\} \in
\mathbb{N}
$. To simplify notation, let us put $\sigma (i,d)=\left(\sum_{k=1}^{i}k^{\lambda }\left\vert
d_{k}-d_{k-1}\right\vert ^{\lambda }\right) ^{\frac{1}{\lambda }}$. Minkowski's inequality and %
\eqref{eq8} yield
\begin{gather*}
\sigma \left(i,d^{(n_{i}^{\ast})}-d\right)=\left( \frac{1}{i}%
\sum_{k=1}^{i}k^{\lambda }\left\vert \left( d_{k}^{(n_{i}^{\ast})}-d_{k}\right) +\left(
d_{k-1}-d_{k-1}^{(n_{i}^{\ast})}\right) \right\vert ^{\lambda }\right) ^{\frac{1}{%
\lambda }} \\
\leq  \left( \frac{1}{i}\sum_{k=1}^{i}k^{\lambda }\left\vert
d_{k}^{(n_{i}^{\ast})}-d_{k}\right\vert ^{\lambda }\right) ^{\frac{1}{\lambda }}+ \left( \frac{1}{i}\sum_{k=1}^{i}k^{\lambda }\left\vert
d_{k-1}-d_{k-1}^{(n_{i}^{\ast})}\right\vert ^{\lambda }\right) ^{\frac{1}{\lambda }} \\
\leq \left( \frac{1}{i}%
\sum_{k=1}^{i}k^{\lambda }\frac{\varepsilon ^{\lambda }}{6^{\lambda }\left(
k+1\right) ^{\lambda }2^{k+1}}\right) ^{\frac{1}{\lambda }}+\left( \frac{1}{i}\sum_{k=1}^{i}k^{\lambda }\frac{\varepsilon
^{\lambda }}{6^{\lambda }k^{\lambda }2^{k}}\right) ^{\frac{1}{\lambda }} \leq \frac{\varepsilon}{3} \text{.}
\end{gather*}%
Taking into account \eqref{eq7}, we get
\begin{gather*}
\sigma \left( i,d^{\left( m\right) }-d\right) \leq \sigma \left( i,d^{\left(
m\right) }-d^{\left( n_{i}^{\ast }\right) }\right) +\sigma \left(
i,d^{\left( n_{i}^{\ast }\right) }-d\right) \\
\leq \left\Vert d^{\left(
m\right) }-d^{\left( n_{i}^{\ast }\right) }\right\Vert _{\mathscr{B}^{\lambda }}+\frac{%
\varepsilon }{3}<\frac{2\varepsilon }{3}\text{ \ \ }\left( \forall m\geq
n_{0}^{\ast \ast }=\max \left\{ n_{0},n_{0}^{\ast }\right\} \text{ \ and }%
\forall i\in
\mathbb{N}
\right)
\end{gather*}
Therefore,
\[
\left\Vert d^{\left( m\right) }-d\right\Vert _{\mathscr{B}^{\lambda }}=\left\Vert
d^{\left( m\right) }-d\right\Vert _{l^{\infty }}+\underset{i}{\sup } \ \sigma
\left( i,d^{\left( m\right) }-d\right) <\frac{\varepsilon }{3}+\frac{%
2\varepsilon }{3}\text{ =}\varepsilon \text{ \ \ }\left(\forall m\geq n_{0}^{\ast
\ast }\right) \text{.}
\]

Finally,
\[
\left\Vert d\right\Vert _{\mathscr{B}^{\lambda }}\leq \left\Vert d-d^{\left( n_{0}^{\ast \ast }\right) }\right\Vert _{\mathscr{B}^{\lambda }}+\left\Vert d^{\left( n_{0}^{\ast \ast }\right) }\right\Vert _{\mathscr{B}^{\lambda }}<\infty \text{.}
\]
Hence, $d \in \mathscr{B}^{\lambda}$.
\end{proof}
\end{thm}

\section{Local to global: behaviour of trigonometric series of a special type}

\begin{defnx}
\label{d1} A sequence of positive numbers $\left\{ d_{n}\right\} $ is said
to be almost decreasing if there exists a constant $M$ such that $%
d_{n+1}\leq Md_{n}$ holds for every $n\in
\mathbb{N}
$. $M$ is the index of almost monotonicity of $\left\{ d_{n}\right\} $. The
space of almost decreasing sequences with index $M$ is denoted by $\mathcal{A%
}_{M}\mathcal{M}$. If $d_{n+1}\leq Md_{n}$ holds true starting from some
integer $n>1$, the corresponding space is denoted by $\mathcal{GA}_{M}\mathcal{M}
$.
\end{defnx}

\begin{rem}
Note that $\mathcal{A}_{1}\mathcal{M}=\mathcal{M}$ is the space of decreasing sequences.
\end{rem}

The role of almost decreasing sequences is nicely illustrated by the following theorem.

\begin{thmx}
Let $\rho _{n}=\sqrt{a_{n}^{2}+b_{n}^{2}}$, $n \in \mathbb{N}$, form an almost decreasing sequence and let
$$\sum A_{n}(x)\equiv \frac{a_0}{2} + \sum a_n \cos nx + b_n \sin nx ,$$
$$\sum B_{n}(x)\equiv \sum a_n \sin nx - b_n \cos nx.$$

\begin{description}[leftmargin=0cm]
\item[(a)] (cf. \cite[Theorem 2]{PI}) If one of the series $\sum A_{n}(x)$ or $\sum B_{n}(x)$ is absolutely convergent at two points $x_0$ and $x_1$ with $|x_0 - x_1| \not\equiv 0 \pmod{\pi}$, then $\sum \rho _{n} < \infty$.
\item[(b)] (cf. \cite[Theorem 2.2]{AV1}) If one of the series $\sum A_{n}(x)$ or $\sum B_{n}(x)$ is $\left[ I\right] _{\lambda }$, $\lambda \geq 1$, convergent at two points $x_0$ and $x_1$ with $|x_0 - x_1| \not\equiv 0 \pmod{\pi}$, then $n \rho _{n}  \rightarrow 0 \ \left[ C_{1}\right] _{\lambda }$. If $\lambda > 1$, then $\sum A_{n}(x)$ is the Fourier series of a function $f\in \underset{1\leq p<\infty }{\bigcap }L^{p}$, $\left[ I\right] _{\lambda }$ convergent to $f$ a.e., and $\sum B_{n}(x)$ is the Fourier series of its conjugate function $\tilde f$, $\left[ I\right] _{\lambda }$ convergent to $\tilde f$ a.e.
\end{description}
\end{thmx}

In the next theorem, we shall replace the condition of $\lambda -$strong convergence by $\lambda -$strong boundedness. A series is said to be $\lambda-$strongly bounded if the sequence of its partial sums is $\lambda-$strongly bounded.

\begin{thm}
\label{T32}
Let $\rho _{n}=\sqrt{a_{n}^{2}+b_{n}^{2}}$, $n \in \mathbb{N}$, form an almost decreasing sequence. If one of the series
$$\sum A_{n}(x)\equiv \frac{a_0}{2} + \sum a_n \cos nx + b_n \sin nx,$$
$$\sum B_{n}(x)\equiv \sum a_n \sin nx - b_n \cos nx$$
is $\lambda -$strongly bounded, $\lambda >1$, at two points $x_{0}$ and $x_{1}$, $\left\vert x_{0}-x_{1}\right\vert
\not\equiv 0\pmod{\pi}$, then $\sum A_{n}(x)$ and $\sum B_{n}(x)$ are Fourier series of functions $f$, $\tilde{f}$, resp., belonging to $L^{p}$ for each $1\leq p < \infty$.

\begin{proof}
Let $\theta _{n}$ be chosen such that $\sin \theta _{n}=\frac{a_{n}}{\rho
_{n}}$ and $\cos \theta _{n}=\frac{b_{n}}{\rho _{n}}$. Then, $\sum A_{n}(x)$ may be written in the form $%
\sum \rho _{n}\sin \left( nx+\theta _{n}\right) $ and $\sum B_{n}(x)=-\sum \rho _{n}\cos \left( nx+\theta _{n}\right)$. Assume that $\sum A_{n}(x)$ is $\lambda-$strongly bounded at two points:
\begin{equation}
\sum_{k=1}^{n}k^{\lambda }\rho _{k}^{\lambda }\left\vert \sin \left(
kx_{i}+\theta _{k}\right) \right\vert ^{\lambda }=O\left( n\right) \text{ as
}n\rightarrow \infty \text{, for }i=0,1\text{.}  \label{eq0}
\end{equation}%
Let $h=x_{0}-x_{1}$. Then $nh=\left( nx_{0}+\theta _{n}\right) -\left(
nx_{1}+\theta _{n}\right) $ and%
\begin{equation*}
\sin nh=\sin \left( nx_{0}+\theta _{n}\right) \cos \left( nx_{1}+\theta
_{n}\right) -\cos \left( nx_{0}+\theta _{n}\right) \sin \left( nx_{1}+\theta
_{n}\right) \text{.}
\end{equation*}%
Therefore,
\begin{equation*}
\left\vert \sin nh\right\vert ^{\lambda }\leq 2^{\lambda }\left( \left\vert
\sin \left( nx_{0}+\theta _{n}\right) \right\vert ^{\lambda }+\left\vert
\sin \left( nx_{1}+\theta _{n}\right) \right\vert ^{\lambda }\right) \text{.}
\end{equation*}%
The last inequality and \eqref{eq0} imply
\begin{equation}
\sum_{k=1}^{n}k^{\lambda }\rho _{k}^{\lambda }\left\vert \sin kh\right\vert
^{\lambda }=O\left( n\right) \text{ as }n\rightarrow \infty \text{.}
\label{eq11}
\end{equation}%
Let $\left\{ \rho _{k}\right\} \in \mathcal{GA}_{M}\mathcal{M} $, $M > 1$. There exists $K\in\mathbb{N}$ such that
\begin{equation*}
\rho _{k-1}\geq \frac{1}{M}\rho _{k}\text{ for }k\geq K\text{.}
\end{equation*}
One has%
\begin{gather*}
\left( k-1\right) \rho _{k-1}\left\vert \sin \left( k-1\right) h\right\vert
+k\rho _{k}\left\vert \sin kh\right\vert \geq \frac{1}{M}\left( k-1\right)
\rho _{k}\left\vert \sin \left( k-1\right) h\right\vert +k\rho
_{k}\left\vert \sin kh\right\vert  \\
\geq \frac{k-1}{Mk}k\rho _{k}\left( \left\vert \sin \left( k-1\right)
h\right\vert +\left\vert \sin kh\right\vert \right)  \\
\geq \frac{1-\epsilon }{M}k\rho _{k}\left( \left\vert \sin \left( k-1\right)
h\right\vert +\left\vert \sin kh\right\vert \right)
\end{gather*}%
for $k>k_{0}= \max \{K, \lfloor\frac{1}{\epsilon }\rfloor \} $, $\epsilon >0$ arbitrarily
small. This and
\begin{gather*}
\left\vert \sin \left( k-1\right) h\right\vert +\left\vert \sin
kh\right\vert \geq \sin ^{2}\left( k-1\right) h+\sin ^{2}kh \\
=1-\cos h\cos \left( 2k-1\right) h\geq 1-\left\vert \cos h\right\vert
\end{gather*}%
yield
\begin{equation*}
\left( k-1\right) \rho _{k-1}\left\vert \sin \left( k-1\right) h\right\vert
+k\rho _{k}\left\vert \sin kh\right\vert \geq M_{1}k\rho _{k}
\end{equation*}%
for $k > k_{0}$, where $M_{1}=M_{1}\left( h\right) >0$. Hence,%
\begin{gather*}
\sum_{k=1}^{n}k^{\lambda }\rho _{k}^{\lambda }=\sum_{k=1}^{k_{0}}k^{\lambda
}\rho _{k}^{\lambda }+\sum_{k=k_{0}+1}^{n}k^{\lambda }\rho _{k}^{\lambda } \\
\leq \sum_{k=1}^{k_{0}}k^{\lambda }\rho _{k}^{\lambda }+\frac{2^{\lambda }}{%
M_{1}^{\lambda }}\sum_{k=k_{0}+1}^{n}\left[\left( k-1\right) ^{\lambda }\rho
_{k-1}^{\lambda }\left\vert \sin \left( k-1\right) h\right\vert ^{\lambda
}+k^{\lambda }\rho _{k}^{\lambda }\left\vert \sin kh\right\vert ^{\lambda }\right]%
\text{.}
\end{gather*}%
Since the first summand in the last line is a finite sum and the second one
is $O\left( n\right) $ as $n\rightarrow \infty $ by \eqref{eq11}, we get
\begin{equation*}
u_{n, \lambda}:=\sum_{k=1}^{n}k^{\lambda }\rho _{k}^{\lambda }=O\left( n\right) \text{
as }n\rightarrow \infty \text{.}
\end{equation*}
Now, let $p>\max \left\{ \frac{\lambda }{\lambda -1},2\right\} $ and $\frac{1%
}{q}=1-\frac{1}{p}$. It is straightforward that $1<q<\min \{\lambda ,2\}$.
Notice that $\sum_{k=1}^{n} k^{\lambda }\rho _{k}^{\lambda }=O\left( n\right) $
implies $u_{n, q}=\sum_{k=1}^{n} k^{q}\rho _{k}^{q}=O\left( n\right) $.

Abel's partial summation formula gives us
\begin{equation}
\label{eq10}
\begin{split}
\sum_{k=1}^{n}\rho _{k}^{q }=\sum_{k=1}^{n}\frac{k^{q }\rho
_{k}^{q }}{k^{q }}=\sum_{k=1}^{n}\frac{u_{k,q}-u_{k-1,q}}{k^{q
}}=\frac{u_{n,q}}{n^{q }}+\sum_{k=1}^{n-1}u_{k,q}\left( \frac{1}{%
k^{q }}-\frac{1}{(k+1)^{q }}\right) \\
=O\left( \frac{1}{n^{q -1}}\right) +O\left( \sum_{k=1}^{n-1}\frac{1}{%
k^{q }}\right) =O\left( 1\right) \text{ as }n\rightarrow \infty \text{.%
}
\end{split}
\end{equation}
By the Hausdorff-Young theorem \cite[(2.3), (ii), p. 101]{ZY}, there exists $f \in L^{p}$ such that $\sum A_{n}(x)$ is the Fourier series of $f$. This and the uniqueness property of Fourier series yield that $f$ belongs to all $L^p$ spaces, $1\leq p < \infty$. Then $\sum B_{n}(x)$ is the Fourier series of $\tilde{f} \in \underset{1\leq p<\infty }{\bigcap }L^{p}$.

The proof is analogous if $\sum B_{n}(x)$ is $\lambda-$strongly bounded at $x_0, x_1$ or if $\sum A_{n}(x)$  is $\lambda-$strongly bounded at $x_0$ and $\sum B_{n}(x)$ at $x_1$.
\end{proof}
\end{thm}

\begin{defnx}
\label{d2} A sequence of positive numbers $\left\{ d_{n}\right\} $ is said
to be quasi decreasing if there exists $\alpha >0$ such that $\left\{
d_{n}/n^{\alpha }\right\} $ is a decreasing sequence starting from some
integer $n\geq 1$. $\alpha $ is the index of quasimonotonicity of $\left\{
d_{n}\right\} $. The space of quasi decreasing sequences with index $\alpha $
is denoted by $\mathcal{Q}_{\alpha }\mathcal{M}$.
\end{defnx}

As an application of the concept introduced by Definition \ref{d2}, we cite the next result.

\begin{thmx} (\cite[Theorem 3.1]{AV1})
Let $\rho _{n}=\sqrt{a_{n}^{2}+b_{n}^{2}}$, $n \in \mathbb{N}$, form a quasi decreasing sequence with index $0<\alpha <1$. Let a trigonometric series $\frac{a_0}{2} + \sum a_n \cos nx + b_n \sin nx $ be strongly convergent at two points $x_0$ and $x_1$ with $|x_0 - x_1| \not\equiv 0 \pmod{\pi}$. Then this series and its conjugate are Fourier series, strongly convergent a.e.
\end{thmx}

Having in mind that the classes of $\lambda-$strongly bounded sequences, $\lambda >1$, are contained in the class of strongly bounded (i.e., $1-$strongly bounded) sequences, we turn a closer attention to the latter case.

We shall consider a new class of logarithmic quasi decreasing sequences.

\begin{defn}
\label{d3} A sequence of positive numbers $\left\{ d_{n}\right\} $ is said
to be logarithmic quasi decreasing if there exists $\beta >0$ such that $%
\left\{ d_{n}/\log ^{\beta }n\right\} $ is a decreasing sequence starting
from some integer $n\geq 2$. $\beta $ is the index of logarithmic
quasimonotonicity of $\left\{ d_{n}\right\} $. The set of logarithmic quasi
decreasing sequences with index $\beta $ is denoted by $\mathcal{L}_{\beta }%
\mathcal{QM}$.
\end{defn}

\begin{thm}\label{T34}
Let $\rho _{n}=\sqrt{a_{n}^{2}+b_{n}^{2}}$, $n\in \mathbb{N}$, form a logarithmic quasi decreasing sequence with index $\beta >1$ ($\rho _{n}\in
\mathcal{L}_{\beta }\mathcal{QM}$). If one of the series
$$\sum A_{n}(x)\equiv \frac{a_{0}}{2}+\sum a_{n}\cos nx+b_{n}\sin nx,$$
$$\sum B_{n}(x)\equiv \sum a_{n}\sin nx-b_{n}\cos nx$$
is strongly bounded at two points $x_{0}$, $x_{1}$, $\left\vert x_{0}-x_{1}\right\vert
\not\equiv 0\pmod{\pi}$, then $\sum \frac{n\rho _{n}^{2}}{%
\log ^{\beta }n}<\infty $. $\sum A_{n}(x)$ and $\sum B_{n}(x)$ are Fourier series of functions $f$, $\tilde{f}$
which belong to $L^{p}$ for each $1\leq p<\infty $.

\begin{proof}
Since $\left\{ \rho _{k}\right\}\in \mathcal{L}_{\beta }\mathcal{QM}$, we have
\begin{equation*}
\rho _{k-1}\geq \frac{\log ^{\beta }\left( k-1\right) }{\log ^{\beta }k}\rho
_{k}\text{ for }k\geq K \geq 2\text{.}
\end{equation*}%
Reasoning as in the proof of Theorem \ref{T32}, we get
\begin{equation*}
u_{n}:=\sum_{k=1}^{n}k\rho _{k}=O\left( n\right) \text{ as }n\rightarrow
\infty \text{.}
\end{equation*}

Now, for any $\alpha >1$, one has
\begin{gather}
\label{eq1}
\sum_{k=2}^{n}\frac{\rho _{k}}{\log ^{\alpha }k}=\sum_{k=2}^{n}\frac{%
k\rho _{k}}{k\log ^{\alpha }k}=\sum_{k=2}^{n}\frac{u_{k}-u_{k-1}}{k\log
^{\alpha }k}\\
=\frac{u_{n}}{n\log ^{\alpha }n}-\frac{u _{1}}{2\log ^{\alpha }2}%
+\sum_{k=2}^{n-1}u_{k}\left( \frac{1}{k \log ^{\alpha
}k }-\frac{1}{(k+1)\log ^{\alpha }(k+1)}\right) \text{.}  \notag
\end{gather}

Obviously
\begin{equation}
\frac{u_{n}}{n\log ^{\alpha }n}=o\left( 1\right) \text{ as }n\rightarrow
\infty \text{.}  \label{eq2}
\end{equation}%
Notice that%
\begin{equation*}
\frac{1}{k \log ^{\alpha }k }-\frac{1}{(k+1)\log
^{\alpha }(k+1)}=\frac{1}{\xi _{k}^{2}\log ^{\alpha }\xi _{k}}\left( 1+\frac{%
\alpha }{\log \xi _{k}}\right) \text{,}
\end{equation*}%
where $\xi _{k}\in \left( k,k+1\right) $. From $\frac{1}{\xi _{k}^{2}\log
^{\alpha }\xi _{k}}\left( 1+\frac{\alpha }{\log \xi _{k}}\right) <\frac{1}{k^{2}\log ^{\alpha }k }\left( 1+\frac{%
\alpha }{\log 2}\right) $ and \\\ $u_{k}=O\left( k \right) $, we get%
\begin{equation*}
u_{k}\left( \frac{1}{k \log ^{\alpha }k }-%
\frac{1}{(k+1)\log ^{\alpha }(k+1)}\right) =O\left( \frac{1}{k \log
^{\alpha }k }\right) \text{.}
\end{equation*}%
Thus,
\begin{equation}
\sum_{k=2}^{n-1}u_{k}\left( \frac{1}{k \log ^{\alpha }k }-\frac{1}{(k+1)\log ^{\alpha }(k+1)}\right) =O\left(\sum_{k=2}^{n-1}
\frac{1}{k\log ^{\alpha }k}\right) =O\left( 1\right) \label{eq3}
\end{equation}
as $n\rightarrow \infty$.\\
\noindent The relations \eqref{eq1}, \eqref{eq2} and \eqref{eq3} yield
\begin{equation}
\sum_{k=2}^{n}\frac{\rho _{k}}{\log ^{\alpha }k}=O\left( 1\right) \text{
as }n\rightarrow \infty \text{, for }\alpha >1\text{.}  \label{eq4}
\end{equation}%
In particular, the series $\sum_{k=2}^{\infty }%
\frac{\rho _{k}}{\log ^{\beta }k}$ is convergent. This and the fact that the
sequence $\left\{ \frac{\rho _{k}}{\log ^{\beta }k}\right\} $ is decreasing
yield $\frac{k\rho _{k}}{\log ^{\beta }k}=o\left( 1\right) $ as $%
k\rightarrow \infty $ by Olivier's theorem.
Now,
\begin{eqnarray*}
\sum_{k=2}^{n}\frac{k\rho _{k}^{2}}{\log ^{\beta }k} &=&\sum_{k=2}^{n}\frac{%
\rho _{k}}{\log ^{\beta }k}\left( u_{k}-u_{k-1}\right) \\
&=&\frac{u_{n}\rho _{n}%
}{\log ^{\beta }n}+\sum_{k=2}^{n-1}u_{k}\left( \frac{\rho _{k}}{\log ^{\beta
}k}-\frac{\rho _{k+1}}{\log ^{\beta }\left( k+1\right) }\right) -\frac{\rho
_{2}u_{1}}{\log ^{\beta }2} \\
&\leq &o\left( 1\right) +C\sum_{k=2}^{n-1}k\left( \frac{\rho _{k}}{\log
^{\beta }k}-\frac{\rho _{k+1}}{\log ^{\beta }\left( k+1\right) }\right)  \\
&=&o\left( 1\right) +C\sum_{k=2}^{n-1}\left( \frac{k\rho _{k}}{\log ^{\beta
}k}-\frac{\left( k+1\right) \rho _{k+1}}{\log ^{\beta }\left( k+1\right) }%
\right) +C\sum_{k=2}^{n-1}\frac{\rho _{k+1}}{\log ^{\beta }\left( k+1\right)
} \\
&=&o\left( 1\right) +C\left( \frac{2\rho _{2}}{\log ^{\beta }2}-\frac{n\rho
_{n}}{\log ^{\beta }n}\right) +C\sum_{k=2}^{n-1}\frac{\rho _{k+1}}{\log
^{\beta }\left( k+1\right) }=O\left( 1\right) \text{ as }n\rightarrow \infty
\text{.}
\end{eqnarray*}%
This proves the first assertion
\begin{equation}
\label{eq12}
\sum_{k=2}^{\infty }\frac{k\rho _{k}^{2}}{\log ^{\beta }k}<\infty \text{.}
\end{equation}%

Concerning the second assertion, \eqref{eq12} and the Riesz-Fischer theorem yield that $\sum A_{n}(x)$ and $\sum B_{n}(x)$ are Fourier series of $f, \tilde{f} \in L^{2}$.

Now, let $p>2$ and $\frac{1}{p}+\frac{1}{q}=1$. Obviously, $1<q<2$. As above,
\begin{equation}
\begin{split}
\overset{n}{\underset{k=1}{\sum }}\rho _{k}^{q}&=\overset{n}{\underset{k=1}{%
\sum }}\frac{1}{k^{q}}k^{q}\rho _{k}^{q}=\overset{n-1}{\underset{k=1}{\sum }}%
\left( \Delta \frac{1}{k^{q}}\right) \overset{k}{\underset{i=1}{\sum }}%
i^{q}\rho _{i}^{q}+\frac{1}{n^{q}}\overset{n}{\underset{i=1}{\sum }}%
i^{q}\rho _{i}^{q}\\
&=O\left( \overset{n-1}{\underset{k=1}{\sum }}%
\frac{1}{k^{q+1}}\overset{k}{\underset{i=1}{\sum }}i^{q}\rho _{i}^{q}\right)+\frac{1%
}{n^{q}}\overset{n}{\underset{i=1}{\sum }}i^{q}\rho _{i}^{q}\text{.}
\end{split}
\label{eq02}
\end{equation}%
Since $\frac{k\rho _{k}}{\log ^{\beta }k}=o\left( 1\right) $ as $%
k\rightarrow \infty $, we have that $k^{q}\rho _{k}^{q}=o\left( \log ^{\beta
q}k\right) $. Therefore,%
\begin{equation*}
\overset{m}{\underset{i=1}{\sum }}i^{q}\rho _{i}^{q}=O\left( m\log ^{\beta
q}m\right) \text{ for }m\in
\mathbb{N}
\text{.}
\end{equation*}%
The last equality, relation \eqref{eq02} and the fact that $q > 1$ yield%
\begin{equation*}
\overset{n}{\sum_{k=1}}\rho _{k}^{q}=O\left(\overset{n-1}{\underset%
{k=1}{\sum }} \frac{\log ^{\beta q}k}{k^{q}}\right) +O\left( \frac{%
\log ^{\beta q}n}{n^{q-1}}\right) =O\left( 1\right) \text{ as }n \rightarrow \infty \text{.}
\end{equation*}%
Thus, $f, \tilde{f} \in \underset{1\leq p<\infty }{\bigcap }L^{p}$ (cf. the end of the proof of Theorem \ref{T32}).
\end{proof}
\end{thm}

\medskip

\begin{rem} Pointwise convergence a.e. of the series $\sum A_{n}(x)$ and $\sum B_{n}(x)$ in Theorem \ref{T34} follows, of course, from the Carleson-Hunt theorem. However, the Kolmogorov-Selyverstov-Plessner theorem \cite[p. 332]{NKB} already serves the purpose since
$$\sum_{k=2}^{\infty}\rho _{k}^{2}\log k < \sum_{k=2}^{\infty}\frac{k\rho _{k}^{2}}{\log ^{\beta }k} < \infty$$
by \eqref{eq12}.
\end{rem}

\medskip
The following remark concerns the relationship between various sequence spaces considered in this paper.

\begin{rem}\label{rem36}
For $0<M_{1}<1<M_{2}$, one has
\begin{equation*}
\mathcal{A}_{M_{1}}\mathcal{M}\subset \mathcal{M}\subset \underset{\beta >0}{%
\cap }\mathcal{L}_{\beta }\mathcal{QM\subset }\underset{\beta >0}{\cup }%
\mathcal{L}_{\beta }\mathcal{QM\subset }\underset{\alpha >0}{\cap }\mathcal{Q%
}_{\alpha }\mathcal{M\subset }\underset{\alpha >0}{\cup }\mathcal{Q}_{\alpha
}\mathcal{M\subset GA}_{M_{2}}\mathcal{M}\text{.}
\end{equation*}

\begin{proof}
It is obvious that $\mathcal{A}_{M_{1}}\mathcal{M}\subset \mathcal{M}\subset
\underset{\beta >0}{\cap }\mathcal{L}_{\beta }\mathcal{QM}$ since $0<M_{1}<1<%
\frac{\log ^{\beta }\left( n+1\right) }{\log ^{\beta }n}$ for any $\beta >0$
and $n\in
\mathbb{N}
$. The inclusions $\underset{\beta >0}{\cap }\mathcal{L}_{\beta }\mathcal{%
QM\subset }\underset{\beta >0}{\cup }\mathcal{L}_{\beta }\mathcal{QM}$ and $%
\underset{\alpha >0}{\cap }\mathcal{Q}_{\alpha }\mathcal{M\subset }\underset{%
\alpha >0}{\cup }\mathcal{Q}_{\alpha }\mathcal{M}$ are trivial. The inclusion $%
\underset{\alpha >0}{\cup }\mathcal{Q}_{\alpha }\mathcal{M\subset GA}_{M_{2}}%
\mathcal{M}$ follows from $\frac{\left( n+1\right) ^{\alpha }}{n^{\alpha }}%
<M_{2}$ for any $\alpha >0$, $M_{2}>1$ and sufficiently large $n\in \mathbb{N}$. Finally, to establish $\underset{\beta >0}{\cup }\mathcal{L}_{\beta }\mathcal{%
QM\subset }\underset{\alpha >0}{\cap }\mathcal{Q}_{\alpha }\mathcal{M}$, it
is enough to check that
\begin{equation}
\frac{\log ^{\beta }\left( n+1\right) }{\log ^{\beta }n}\leq \frac{\left(
n+1\right) ^{\alpha }}{n^{\alpha }}  \label{eq5}
\end{equation}%
holds true for $\beta >0$, $\alpha >0$ and $n$ sufficiently large. The last
inequality is equivalent to
\begin{equation*}
\frac{\log \left( n+1\right) }{\log n}\leq \left( 1+\frac{1}{n}\right)
^{\gamma }
\end{equation*}%
where we put $\gamma =\frac{\alpha }{\beta }>0$. Subtracting $1$ from both
sides, we get%
\begin{equation*}
\frac{\log \left( 1+\frac{1}{n}\right) }{\log n}\leq \left( 1+\frac{1}{n}%
\right) ^{\gamma }-1\text{.}
\end{equation*}%
According to Taylor's formula, the left hand side is equal to $\frac{1}{n\log n}-\frac{1%
}{2n^{2}\log n}+O(\frac{1}{n^{3}\log n})$, while the right hand side is equal to
$\frac{\gamma }{n}+\frac{\gamma \left( \gamma -1\right) }{2n^{2}}+O\left(
\frac{1}{n^{3}}\right)$. Therefore, inequality \eqref{eq5} holds true for $%
\beta >0$, $\alpha >0$ and sufficiently large $n$.
\end{proof}
\end{rem}

\begin{rem}
In Remark \ref{rem36} we are actually dealing with equivalence classes. Namely, while proving the inclusions, we suppose
that $\left\{ d_{k}\right\} _{k\geq k_{0}}$ and $\left\{
d_{k}\right\} _{k\geq k_{1}}$, $k_{0}\neq k_{1}$, represent the same sequence.
\end{rem}

\begin{rem}
We have seen in Theorem \ref{T32} that if $\left\{ \rho _{k}\right\} \in \mathcal{GA}_{M}\mathcal{M} $, $M > 1$, then a mere $\lambda-$boundedness, $\lambda >1$, of the series $\sum A_{n}(x)$ or $\sum B_{n}(x)$ at two distinct points is sufficient to conclude that these series are Fourier series of functions $f, \tilde{f}$ belonging to all $L^{p}$ spaces, $1\leq p < \infty$. In the case $\lambda =1$, the same conclusion is valid under a stronger assumption $\rho _{k}\in \mathcal{L}_{\beta }\mathcal{QM}$, $\beta > 1$. For intermediate classes $\mathcal{Q}_{\alpha }\mathcal{M}$, $\alpha > 0$, the same techniques of the proof yield the following theorem.
\end{rem}

\begin{thm}
Let $\left\{ \rho _{k}\right\} \in \mathcal{Q}_{\alpha }\mathcal{M}$, $\alpha > 0$. If $\sum A_{n}(x)$ or $\sum B_{n}(x)$ is strongly bounded at two points $x_{0}$ and $x_{1}$, $\left\vert x_{0}-x_{1}\right\vert
\not\equiv 0\pmod{\pi}$, then $\sum k^{1-\alpha}\rho _{k}^{2} < \infty$. These series are Fourier series of functions $f$, $\tilde{f}$ which belong to $L^{2}$ if $\alpha\in (0,1)$. Moreover, $f, \tilde{f}\in L^{p}$, $2<p<\frac{1}{\alpha}$, if $\alpha\in (0,\frac{1}{2})$.
\end{thm}

\section{Strong convergence and generalized variation}



Given a trigonometric series
\begin{equation}
\frac{a_{0}}{2}+\sum_{k=1}^{n}a_{k}\cos kx+b_{k}\sin kx  \label{eq00}
\end{equation}%
let $s_{n}\left( x\right) $ and $\sigma_{n}\left( x\right) $ denote the ordinary $n-$th partial sum and $n-$th Ces\`{a}ro $(C,1)$ partial sum of \eqref{eq00},
respectively. If \eqref{eq00} is a Fourier series of $f\in L^{1}$, we shall write $s_{n}f$ and $\sigma_{n}f$ for the partial sums $s_{n}$ and $\sigma_{n}$.

We will consider the following classes of functions
$$\mathcal{S}^{\lambda}=\left\{f\in L^{1} : s_{n}f \rightarrow f \ \left[ I\right]_{\lambda } \text{ a.e.}\right\} \text{,}$$
$$\mathscr{S}^{\lambda }=\left\{ f\in C:s_{n}f\rightarrow f \ \left[ I\right]_{\lambda }\text{ uniformly}\right\} \text{,}$$
$$\mathscr{U}=\left\{ f\in C:s_{n}f\rightarrow f \ \text{ uniformly}\right\} \text{,}$$
where $C$ is the space of $2\pi -$periodic continuous functions.

For $\lambda\geq1$, it is known (see \cite{NTM3}) that
$$\mathcal{S}^{\lambda}=\left\{f\in L^{1} : \sum_{k=1}^{n}k^{\lambda}\rho_{k}^{\lambda}=o(n)\right\}$$
and
$$\mathscr{S}^{\lambda }=\left\{ f\in C : \sum_{k=1}^{n}k^{\lambda}\rho_{k}^{\lambda}=o(n)\right\} \text{.}$$

By $W$ we denote the class of \textit{regulated functions}, i.e. functions possessing the one-sided limits at each point. Every regulated function is bounded and has at most a countable set of discontinuities. Regulated functions have a particular role in the matter of everywhere convergence of Fourier series.

Important subclasses of the class $W$ stem from various concepts of generalized bounded variation. In the sequel, let $f(I):=f(b)-f(a)$ for arbitrary subinterval $(a,b)$ of $(0,2\pi)$ and the supremum in defining sums below is always taken over all finite collections of nonoverlapping subintervals $I_i$ of $(0,2\pi)$.

According to N. Wiener \cite{WI}, a function $f$ is of $p-$\textit{bounded variation}, $p\geq 1$, on $\left[ 0,2\pi \right] $ and belongs to the class $V_{p}$ if
\begin{equation*}
V_{p}(f)=\sup \left\{ \underset{i}{\sum }\left\vert f(I_{i})\right\vert
^{p}\right\} ^{1/p}<\infty \text{.}
\end{equation*}

A function $f$ is of $\phi-$\textit{bounded variation} (L. C. Young \cite{YO}) on $\left[ 0,2\pi \right] $ and belongs to the class $V_{\phi }$ if
 \begin{equation*}
V_{\phi }(f)=\sup \left\{ \underset{i}{\sum }\phi \left( \left\vert
f(I_{i})\right\vert \right) \right\} <\infty \text{.}
\end{equation*}
Here, $\phi $ is a continuous function defined on $\left[ 0,\infty \right) $ and strictly
increasing from $0$ to $\infty$.

Notice that by taking $\phi \left( u\right) =u$ we get Jordan's class $BV$, while $\phi
\left( u\right) =u^{p}$ gives Wiener's class $V_{p}$.

A function $f$ is of $\Lambda-$\textit{bounded variation} (D. Waterman \cite{WA}) on $\left[ 0,2\pi \right] $ and belongs to the class $\Lambda BV$ if
\begin{equation*}
V_{\Lambda }(f)=\sup \left\{ \underset{i}{\sum }\left\vert
f(I_{i})\right\vert /\lambda _{i}\right\} <\infty \text{,}
\end{equation*}
where $\Lambda =\left\{ \lambda _{n}\right\} $ is a nondecreasing sequence of positive numbers tending to infinity, such that $\sum 1/\lambda _{n}$ diverges.

In the case when $\Lambda =\left\{ n\right\} $, the sequence of positive integers, the function $f$ is said to be of \textit{harmonic bounded variation} and the corresponding class is denoted by $HBV$.

$BV$ is the intersection of all $\Lambda BV$ spaces and $W$ is the union of all $\Lambda BV$ spaces \cite{PE}.

D. Waterman also introduced the notion of continuity in $\Lambda-$variation  to
provide a sufficient condition for $(C, \alpha)-$summability of Fourier series \cite{WA1}. Let $\Lambda^{m}=\{\lambda_{n+m}\}$, $m=0,1,2,\ldots$. A function $f\in \Lambda BV$ is said to be \textit{continuous in $\Lambda-$variation} (or to belong to $\Lambda_{c}BV$) if $V_{\Lambda^{m}}(f)\rightarrow 0$ as $m\rightarrow \infty$.

Clearly, $\Lambda_{c}BV\subseteq \Lambda BV$. Functions from $\Lambda_{c}BV$ admit much better estimates of their Fourier coefficients (see \cite{WN,SA}).

The \textit{modulus of variation} (Z. Chanturiya \cite{CH})
of a bounded function $f$ is the function $\nu _{f}$ whose domain is the set
of positive integers, given by
\begin{equation*}
\nu _{f}\left( n\right) =\sup \left\{ \underset{k=1}{%
\overset{n}{\sum }}\left\vert f(I_{k})\right\vert \right\} \text{.}
\end{equation*}%
The modulus of variation of any bounded function is nondecreasing and
concave. Given a function $\nu $ whose domain is the set of positive
integers with such properties, then by $V\left[ \nu \right] $ one denotes
the class of functions $f$ for which $\nu _{f}\left( n\right) =O\left( \nu
\left( n\right) \right) $ as $n\rightarrow \infty $. We note that $V_{\phi
}\subseteq V\left[ n\phi ^{-1}\left( 1/n\right) \right] $ and $W=\left\{
f:\nu _{f}\left( n\right) =o\left( n\right) \right\} $ \cite{CH}.

The relationship between Waterman's and Chanturiya's concepts was established in \cite{AV}. M. Avdispahi\'{c} proved the following inclusions between Wiener's, Waterman's
and Chanturiya's classes of functions of generalized bounded variation.

\begin{thmx}[cf. Theorem 4.4. in \cite{AV0}]
\label{AVth}
\begin{equation*}
\left\{ n^{\alpha }\right\} BV\subset V_{\frac{1}{1-\alpha }}\subset V\left[
n^{\alpha }\right] \subset \left\{ n^{\beta }\right\} BV\text{,}
\end{equation*}%
for $0<\alpha <\beta <1$.
\end{thmx}

\medskip
The next two theorems are related to strong convergence and strong boundedness of Fourier series of regulated functions. As always, by $\tilde{f}$ we denote the conjugate function of a function $f$.

\begin{thm} Let $\lambda\geq1$. Then
\label{lema}
\leavevmode
\begin{description}
\item[i)] $W\cap \mathcal{S}^{\lambda }=\mathscr{S}^{\lambda }$.

\item[ii)] If $f,\tilde{f}\in W$, then $f,\tilde{f}\in C$.

\item[iii)] If $f\in \mathcal{S}^{\lambda }$ and $\tilde{f}\in W$, then $\tilde{f}\in %
\mathscr{S}^{\lambda }$.

\item[iv)] If $f\in HBV$ and $\tilde{f}\in W$, then $f,\tilde{f}\in \mathscr{U}$.
\end{description}

\begin{proof}
\leavevmode
i) Let $f$ be an arbitrary function in $W\cap \mathcal{S}^{\lambda }$. Recall that $\mathcal{S}^{\lambda }\subset \mathcal{S}$ \cite[Theorem 1. (iii)]{NTM3}. Thus, $\sum_{k=1}^{n}k \rho _{k}=o\left( n\right)$, as $n \rightarrow \infty$. By \cite[Theorem 3, p. 183 and Corrolary 2, p. 185]{NKB}, $f$ can not have discontinuities of the first kind. It follows that $f$ is a continuous function. Its Fourier series is $(C,1)$ uniformly summable. Therefore, $f\in \mathscr{S}^{\lambda }$. The converse, $\mathscr{S}^{\lambda } \subseteq W\cap \mathcal{S}^{\lambda }$, is trivial. \\
ii) Let $f,\tilde{f}\in W$. If there exists a point $x_{0}$
such that, e.g., $f\left( x_{0}+0\right) -f\left( x_{0}-0\right) >0$, then by \cite[Teorem 8.13, vol. I, p. 60]{ZY} $\tilde{S}_{n}\left(x_{0},f\right) \rightarrow -\infty $. Hence, $\tilde{\sigma}%
_{n}\left( x_{0},f\right) \rightarrow -\infty $, which contradicts the fact that $\tilde{\sigma}_{n}\left( x_{0},f\right) =\sigma
_{n}\left( x_{0},\tilde{f}\right) \rightarrow \frac{1}{2}\left[ \tilde{f}%
\left( x_{0}+0\right) +\tilde{f}\left( x_{0}-0\right) \right] $ \cite[Fej\'{e}r's theorem 3.4, vol. I, p. 89]{ZY}. Therefore, function $f$ is continuous. Analogously, the function $\tilde{f}$ is continuous. \\
iii) Let $\tilde{f}\in W$. The conjugate series is $(C,1)$ summable to $\tilde{f}$ a.e. \cite[p. 524]{NKB}. Therefore, $f\in \mathcal{S}^{\lambda }$ implies $\tilde{f}\in \mathcal{S}^{\lambda }$. Hence, $\tilde{f} \in W\cap \mathcal{S}^{\lambda }=\mathscr{S}^{\lambda }$ by i).\\
iv) By ii) above, $f,\tilde{f}\in C$. Now, $f \in HBV \cap C$ implies uniform convergence of its Fourier series \cite{WA}. However, $\tilde{f}$ being also continuous, its Fourier series is necessarily uniformly convergent as well, by \cite[Theorem 1, p. 592]{NKB}.
\end{proof}
\end{thm}


\begin{thm}
\label{nbv}
\leavevmode
\begin{description}
\item[i)] $\left\{ n^{1/2}\right\} BV\cap C\subset \mathscr{S}^{2}$.
\item[ii)] If $f\in \left\{ n^{1/2}\right\} BV$ and $\tilde{f} \in W$, then $f,\tilde{f}%
\in \mathscr{S}^{2}$.
\item[iii)]If $f\in V_{2}$, then sequence $\left\{ s_{n}f\right\} $ is $2-$strongly bounded.
\end{description}
\begin{proof}
i) Let $f\in \left\{ n^{1/2}\right\} BV\cap C$. Uniform convergence of the Fourier series follows from \cite{WA}. M. Avdispahi\'{c} \cite[Theorem 11.1]{AV0} proved that the condition
\begin{equation}
\frac{1}{n}\sum_{k=1}^{n}k^{2}\rho _{k}^{2}=o\left( 1\right) \text{ as }%
n\rightarrow \infty  \label{eq01}
\end{equation}%
is necessary and sufficient for continuity of $\ f\in \left\{
n^{1/2}\right\} _{c}BV$.\ According to \cite[Theorem 3.1]{FPW}
the equality $\Lambda _{c}BV=\Lambda BV$ holds if and only if $S_{\lambda
}<2$, where $S_{\lambda }$ is the Shao-Sablin index defined by%
\begin{equation*}
S_{\lambda }:=\underset{n\longrightarrow \infty }{\lim \sup }\frac{%
\sum\nolimits_{i=1}^{2n}\frac{1}{\lambda _{i}}}{\sum\nolimits_{i=1}^{n}\frac{%
1}{\lambda _{i}}}
\end{equation*}%
for every proper $\Lambda -$sequence $\Lambda =\left\{ \lambda _{i}\right\} $. In case of $\Lambda =\left\{ i^{1/2}\right\} $, we have%
\begin{equation*}
S_{\lambda }=\underset{n\longrightarrow \infty }{\lim \sup }\frac{%
\sum\nolimits_{i=1}^{2n}\frac{1}{\sqrt{i}}}{\sum\nolimits_{i=1}^{n}\frac{1}{%
\sqrt{i}}}=\underset{n\rightarrow \infty }{\lim }\frac{\int\nolimits_{1}^{2n}%
\frac{dx}{\sqrt{x}}}{\int\nolimits_{1}^{n}\frac{dx}{\sqrt{x}}}=\underset{%
n\rightarrow \infty }{\lim }\frac{\sqrt{2n}-1}{\sqrt{n}-1}=\sqrt{2}<2\text{.}
\end{equation*}%
Therefore, \eqref{eq01} holds for $f\in \left\{ n^{1/2}\right\} BV\cap C$.
Since%
\begin{equation*}
\frac{1}{n}\sum_{k=1}^{n}k^{2}\left\vert s_{k}f-s_{k-1}f\right\vert ^{2}=%
\frac{1}{n}\sum_{k=1}^{n}k^{2}\rho _{k}^{2}\left\vert \sin \left( kx+\theta
_{k}\right) \right\vert ^{2}\leq \frac{1}{n}\sum_{k=1}^{n}k^{2}\rho _{k}^{2}%
\text{,}
\end{equation*}%
\eqref{eq01} and uniform convergence of $\left\{ s_{n}f\right\} $ imply that $\left\{ s_{n}f\right\} $ is $2-$strongly
convergent uniformly, i.e. $f\in \mathscr{S}^{2}$.\\
ii) If $f\in \left\{ n^{1/2}\right\} BV$ and $\tilde{f} \in W$, then $f,\tilde{f}\in C$ by Theorem \ref{lema} ii). Now, $f \in \mathscr{S}^{2}$ according to i) above. Moreover, $\tilde{f}\in \mathscr{S}^{2}$ by Theorem \ref{lema} iii).\\
iii) If $f\in V_{2}$, then $\frac{1}{n}\sum_{k=1}^{n}k^{2}\rho _{k}^{2}=O\left(
1\right) $ \cite[proof of Lema 3.1]{AVS}, and the sequence $\left\{
s_{n}f\right\} $ is $2-$strongly bounded.
\end{proof}
\end{thm}

\begin{rem}
 In view of Theorem \ref{AVth}, the analogues of Theorem \ref{nbv} i) and ii) are valid for Wiener classes $V_{p}$, $1\leq p<2$, and Chanturiya classes $V[n^{\alpha}]$, $0<\alpha<\frac{1}{2}$.
\end{rem}

\end{document}